\magnification=\magstep{1}
%\pagewidth{5.65in}
%\pageheight{7.5in}
\hoffset0.4 true in
\voffset-30pt
%\TagsOnRight
\input amstex
\documentstyle{amsppt}
\nologo
\loadbold
\topmatter
\title {A strong uniqueness theorem \\for planar vector fields}
\endtitle
\author {S. Berhanu and J.Hounie}
\endauthor
\address
Department of Mathematics, Temple University,
Philadelphia, PA 19122-6094, USA
\endaddress
\email
berhanu\@math.temple.edu
\endemail
\address
Departamento de Matem\'atica, UFSCar,
13.565-905,  S\~ao Carlos, SP,  BRASIL
\endaddress
\email
hounie\@ufscar.dm.br
\endemail
\dedicatory
Dedicated to Prof.~Constantine Dafermos on the occasion of his 60th anniversary
\enddedicatory
\thanks
 Work supported in part by CNPq, FINEP, FAPESP and a Research Incentive Fund grant of Temple University
\endthanks
\keywords
Weak boundary values, locally integrable vector fields, uniqueness in the Cauchy problem
\endkeywords
\subjclass Primary 35F15, 35B30, 42B30; Secondary 42A38, 30E25
\endsubjclass

\def\ce{\Bbb C}
\def\erre{\Bbb R}
\def\ccinf{C^\infty_{c}}
\def\cinf{C^\infty}

\define\<{\langle}

\define\>{\rangle}

\redefine \D{\Cal{D}}

\redefine \|{\Vert}
\def\re{\hbox{\rm Re}\, }
\def\im{\hbox{\rm Im}\, }

\abstract{ This work establishes a strong uniqueness property for
a class of planar locally integrable vector fields. A result on
pointwise convergence to the boundary value is also proved for
bounded solutions. }
\endabstract
\endtopmatter
\NoBlackBoxes
 \document
 \heading
 {Introduction}
 \endheading
Consider a  complex, smooth vector field
 $$
L=\frac{\partial}{\partial y}+a(x,y)\frac{\partial}{\partial x}
$$
defined in a neighborhood of the origin in $\erre^2$. We are interested in
the following uniqueness question: if a function $u(x,y)$ defined in a neighborhood of the origin satisfies
$$
\cases
Lu=0&\text{for $y>0$ and}\\
u(x,0)=0,&
\endcases \tag 0.1
$$ can we conclude that $u(x,y)$ vanishes identically in a
neighborhood of the origin? In 1960 P. Cohen [C] (see also [Z] and
the references therein) constructed smooth functions $u(x,y)$ and
$a(x,y)$ defined on the plane such that \roster
\item $Lu(x,y)=(u_y+au_x)(x,y)=0$;
\item $u(x,y)=a(x,y)=0$ for all $y\le0$;
\item $supp\, u=supp\,a=\{(x,y):\quad y\ge0\}$.
\endroster
In particular, $u$ satisfies (0.1) and does not vanish identically
in any neighborhood of the origin which shows that some additional
hypothesis must be made on $L$ if one  wants to obtain uniqueness.
A quite satisfactory answer is known for the class of locally
integrable vector fields ([T1]):  if $u(x,y)=0$ satisfies (0.1)
and $L$ is locally integrable,  $u$ must vanish on a small
rectangle $(-\delta,\delta)\times(0,\delta)$.

In this article we investigate a stronger uniqueness property for
locally integrable vector fields, replacing the condition that
$u(x,0)$ vanish identically by the weaker hypothesis that the
integral of $\ln|u(x,0)|$ be equal to $-\infty$, and consider
one-sided solutions, i.e., $u(x,y)$ is only assumed to satisfy the
equation on one side of the initial curve $\{y=0\}$, conditions
that are classically known to guarantee uniqueness for the
Cauchy-Riemann operator ([F], [RR], see also [Du] and the
references therein). After an appropriate local change of
variables that preserves the initial curve $\{y=0\}$, any elliptic
vector field can be transformed into a multiple of the
Cauchy-Riemann operator and this shows that elliptic vector fields
share this strong uniqueness property. However, this condition is
not enough to  guarantee uniqueness for the vector field
$\partial_y$ so an additional hypothesis has to be made on $L$ if
it is to possess the strong uniqueness property under scrutiny. It
turns out that a much weaker assumption than ellipticity is enough
to ensure that $L$ will share with the Cauchy-Riemann operator
this strong uniqueness property for bounded solutions. All we need to assume is that the integral curve of
$X=\re L$ through the origin contains a sequence of points  on
which $L$ is elliptic and the sequence converges to the origin (see Theorem
1.2 below for the precise statement). It is also shown that this
geometric condition is necessary for the validity of the strong
uniqueness property. The work [J] and the recent article [Co]
contain results on this kind of uniqueness property. In [Co] the
author established the strong uniqueness property for approximate
solutions of a class of planar vector fields. A continuous
function $u$ is said to be an {\it approximate solution } for a
vector field $L$ if $Lu$ is locally integrable and satisfies the
inequality $|Lu|\leq M|u|$ for some constant $M$. In [J] the setup
involves a generic CR manifold $\Cal M$ in $\Bbb C^N$ and a
submanifold $E$ of $\Cal M$ satisfying $T\Cal M + J(T\Cal M)=TE+J(TE)$
over the points of $E$ and where $J$ is the complex structure map.
It was proved that if $u$ is a continuous CR function in a neighborhood of $p\in E$ and if for
$z\in E$,  $|u(z)|\leq h(|z-p|)$ for some continuous,
increasing $h$ on $[0,\infty)$ with $\int_0^1 \ln h(r)\,dr=-\infty$, then $u$
vanishes on the Sussmann orbit through $p$. We emphasize that in
the result of [J], the function $u$ is a solution in a full
neighborhood of  $E$ while in our result (Theorem
1.2), we consider a solution defined only on one side.
\newline This paper is organized as follows. In sections 1 and 2
the proof of the main result Theorem 1.2 is presented. Section 3
is devoted to a theorem on pointwise convergence to weak boundary
values.

\smallskip
\heading {1. A uniqueness criterion for locally integrable vector
fields}
\endheading
We recall that a  vector field with smooth complex coefficients
 $$
L=\frac{\partial}{\partial y}+a(x,y)\frac{\partial}{\partial x}
$$
defined in a neighborhood of the origin in $\erre^2$ is said to be locally integrable at the origin if there exists a smooth function $Z(x,y)$ defined in a neighborhood of the origin such that
\roster
\item $LZ=0$; and
\item $dZ(0,0)\not=0$.
\endroster
A function $Z(x,y)$ satisfying these properties is called a first
integral of $L$. Throughout this paper we assume that $L$ is
locally integrable at the origin. Let $u(x,t)$ be a bounded
measurable function that satisfies $Lu=0$ on $(-a,a)\times(-b,b)$
in the sense of distributions. For a general $L^\infty$ function
the restriction to a horizontal line $t=\text{constant}$ is not
defined because lines have measure zero with respect to the
Lebesgue two dimensional measure. On the other hand, Fubini's
theorem implies that there exist a null set $E\subset(-b,b)$ such
that if $t_0\notin E$ the function $(-a,a)\owns x\mapsto u(x,t_0)$
is measurable. Since $u$ is assumed to satisfy the equation $Lu=0$
we may assert more: since the wave front set $WF(u)$ is contained
in the characteristic set of $L$ $$
\left\{(x,y;\xi,\eta)\in(-a,a)\times(-b,b)\times\erre^2: \quad \im
a(x,y)\xi=0, \,\,\eta +\re a(x,y)\xi=0\right\} $$ it follows
([Hor,Corollary 8.2.7]) that {\it for all } $|t_0|<b$ there is a
well defined restriction of $u$ to the horizontal lines
$t=t_0\in(-b,b)$
---called the trace of $u$ at $t=t_0$--- that will be denoted by
$u(x,t_0)\in\D'(-a,a)$. Furthermore, the distribution $u(x,t_0)$
depends continuously on $t_0$, i.e., if $\phi(x)\in\ccinf((-a,a))$
the function $(-b,b)\owns t\mapsto \<u(x,t),\phi(x)\>$ is
continuous, in fact, it is smooth. We are here using the same
notation, namely $u(x,t_0)$, to denote two different notions: the
distribution trace of $u$ and the pointwise restriction; however
this may not cause confusion as they coincide for almost every $t$
(see for instance [HT,Lemma B.2]).

The most basic and general uniqueness theorem for $L$ is a
by-product of the Baouendi-Treves approximation formula ([BT],
[T3],[T1]):  if $u(x,0)=0$ then $u$ must vanish in a neighborhood
of the origin. This uniqueness result can be strengthened as
follows for one-sided solutions ({\it cf.}, e.g., [T2],
[HM,p.1313]): if $u(x,t)$ is only defined on $(-a,a)\times(0,b)$
where it satisfies $Lu=0$ and we assume that $u(x,0)=0$ in the
sense that $$
\lim_{t\searrow0}\<u(x,t),\phi(x)\>=0,\quad\phi\in\ccinf((-a,a)),
\tag1.1 $$ then $u(x,t)$ must vanish identically for $|x|<\delta$,
$0<t<\delta$ if $\delta>0$ is sufficiently small. For a general
distribution that satisfies $Lu=0$ on $(-a,a)\times(0,b)$ the
limit on the right hand side of (1.1) may not exist but here it
does because $u$ is bounded. Indeed \proclaim{Lemma 1.1} Let $$
L=\frac{\partial}{\partial y}+a(x,y)\frac{\partial}{\partial x},
$$ $a(x,t)\in C^{\infty}$ on $(-a,a)\times(-b,b)$, be a not
necessarily locally integrable vector field. Let $f$ be a bounded
function on $(-a,a)\times(0,b)$ such that $Lf=0$. Then, as
$y\searrow 0$, $f(x,y)$ converges in ${\D}'((-a,a))$ to a bounded
function $bf(x)\doteq\lim_{y\searrow 0}f(x,y)\in L^\infty(-a,a)$.
\endproclaim
\proclaim{Remark}
{\rm This lemma is a variation of Lemma 1.2 in [BH1].
\endproclaim
We postpone the proof of the Lemma ---to be found in the Appendix--- and continue our  discussion of uniqueness. According to our previous discussion, it is known that if $u(x,y)$ is defined and bounded for $y>0$, satisfies the equation $Lu=0$ and its boundary value $bu(x)$ is zero, then $u$ must vanish identically for $|x|<\delta$, $0<t<\delta$, if $\delta$ is small. On the other hand, if $L$ is the Cauchy Riemann operator
$$
\frac{\partial}{\partial \overline z}=\frac{1}{2}\left(\frac{\partial}{\partial x}+i
\frac{\partial}{\partial y}\right)
$$
so $U(x+iy)=u(x,y)$ is bounded and holomorphic for $y>0$ it is well known that it is enough to assume that
$$
\int_{-a}^a\ln|bu(x)|=-\infty \tag 1.2
$$
to conclude that $u$ must vanish. Note that if
\roster
\item "i)" $bu(x)$ vanishes on a set of positive measure, or
\item "ii)" $bu(x)$ has a zero of exponential order, i.e., $|bu(x)|\le A\exp(-B|x-x_0|^{-1})$ for some $x_0\in(-\delta,\delta)$,
\endroster
then (1.2) holds. We wish to extend this finer type of uniqueness
property to a class of locally integrable vector fields. Of
course, the class cannot contain the vector field $L=\partial_y$
which obviously fails to have this type of uniqueness as any
function $u(x)$ independent of $y$ satisfies the homogeneous
equation $\partial_yu=0$. This strong divergence in  behavior
between the Cauchy-Riemann operator and $\partial_y$ is explained
by the fact that the former is elliptic at every point while the
second is not elliptic at any point. Loosely speaking, our class
must exhibit some degree of ellipticity in order to have a chance
to share the uniqueness property we are interested in with the
Cauchy-Riemann operator. A precise geometric property that
characterizes vector fields  $L$ possessing this type of
uniqueness is given by the theorem below. Let's write $L=X+iY$
with $X$ and $Y$ real and note that $L$ is elliptic precisely at
the points where the real vector fields $X$ and $Y$ are linearly
independent. \proclaim{Theorem 1.2} Let $$
L=\frac{\partial}{\partial y}+a(x,y)\frac{\partial}{\partial
x}=X+iY, $$ $a(x,t)\in C^{\infty}$ on $(-a,a)\times(-b,b)$, be
locally integrable. Assume that on the integral curve of $X$ that
passes through the origin there is a sequence of points
$p_n=(x_n,y_n)$ such that \roster
\item $L$ is elliptic at $p_n$, i.e., $X(p_n)$ and $Y(p_n)$ are linearly independent;
\item $y_n>0$ and $p_n\to(0,0).$
\endroster
Then there exists $0<\delta<a$ such that every function $u(x,y)\in
L^\infty((-a,a)\times(0,b))$ that satisfies \roster
\item "(3)" $Lu(x,y)=0$ for $y>0$;
\vskip 3pt
\item "(4)"
$
\displaystyle \int_{-\delta}^\delta \ln|bu(x)|\,dx=-\infty;
$
\endroster
must vanish identically on $(-\delta,\delta)\times(0,\delta)$.

Conversely, if no sequence $p_n=(x_n,y_n)$ on the integral curve
of $X$ that passes through the origin satisfies {\rm (1)} and {\rm
(2)}, there exists a function $u(x,y)\in
C^\infty((-\delta,\delta)\times[0,\delta))$ that satisfies {\rm
(3)} and {\rm (4)} but does not vanish identically on the
intersection of $(-\delta,\delta)\times(0,\delta))$ with any
neighborhood of the origin.
\endproclaim
\demo{Proof} We may find new local coordinates in a neighborhood
of the origin that preserve the $x$-axis and the upper half plane
$\{y>0\}$ in which a first integral $Z(x,y)$ of $L$ has the form
$$ Z(x,y)=x+i\varphi(x,y),\qquad\varphi(0,0)=\varphi_x(0,0)=0, $$
with $\varphi(x,y)$ real. Since $LZ=0$, modulo a nonvanishing
factor, the expression of $L$ in the new coordinates is $$
L=X+iY=\frac{\partial}{\partial
y}-\frac{i\varphi_y}{1+i\varphi_x}\frac{\partial}{\partial
x},\qquad $$ $$ X=\frac{\partial}{\partial y}-\frac{\varphi_y
\varphi_x}{1+\varphi_x^2}\frac{\partial}{\partial x}, \qquad
Y=-\frac{\varphi_y }{1+\varphi_x^2}\frac{\partial}{\partial x}. $$
We see that $X$ and $Y$ become linearly independent precisely at
the points where $\varphi_y$ vanishes. Note that
$y\mapsto\varphi(0,y)$ cannot vanish identically on any interval
$[0,\epsilon]$, $0<\epsilon<b$. Indeed, if it did we would
conclude that $X=\partial_y$ on the segment
$\{0\}\times[0,\epsilon)$ which would then be an integral curve of
$X$ on which $Y$ vanishes, contradicting hypothesis (1). For
$0<\delta<b$ to be determined later (we will have to shrink
$\delta$ several times) set $$ M(x)=\sup_{0\le y\le
\delta}\varphi(x,y),\qquad m(x)=\inf_{0\le y\le
\delta}\varphi(x,y). $$ Then $M(x)$ and $m(x)$ are Lipschitz
continuous functions, $m(x)\le M(x)$ and $m(0)<M(0)$. Assume that
$u(x,y)\in L^\infty((-a,a)\times(-0,b))$  satisfies hypotheses (3)
and (4) and we wish to show that $u$ vanishes near the origin. Now
the Baouendi-Treves approximation formula ([BT], [HM,Cor2.2])
furnishes a sequence of holomorphic polynomials $P_k(\xi+i\eta)$
such that the functions $u_k(x,y)=P_k(Z(x,y))$ satisfy the
following properties for some fixed $\delta>0$ and $K>0$: \roster
\item $|u_k(x,y)|\le K$ on $(-\delta,\delta)\times(0,\delta)$;
\item $u_k(x,y)\to u(x,y)$ a.e. on $(-\delta,\delta)\times(0,\delta)$;
\item $u_k(x,0)\to bu(x)$ a.e. on $(-\delta,\delta)$.
\endroster
Let's denote by $Q$ the rectangle
$(-\delta,\delta)\times(0,\delta)$ and by $\Omega$ the region of
the complex plane bounded by the Lipschitz curves $M(x)$ and
$m(x)$, more precisely, $$
\Omega=\left\{\zeta=\xi+i\eta:\quad|\xi|<\delta,\,\,
m(\xi)<\eta<M(\xi) \right\}, $$ where we have shrunk $\delta>0$ to
ensure that $m(x)<M(x)$ for $|x|<\delta$. It follows from (1) that
the sequence $P_k(\zeta)$ is uniformly bounded on $\Omega$ and by
Montel's theorem we may assume, after passing to a subsequence,
that $P_k(\zeta)$ converges uniformly over compact subsets of
$\Omega$ to a bounded holomorphic function $U(\zeta)$. Since the
functions $u_k(x,t)=P_k\circ Z(x,t)$, $k=1,2\dots,$ satisfy
$Lu_k=0$, converge pointwise  to $u^*(x,t)=U\circ Z(x,t)$ on
$Z^{-1}(\Omega)$ and the $u_k$ are uniformly bounded, we see that
$u^*(x,t)$ also satisfies $Lu^*=0$ on $Z^{-1}(\Omega)$.
Furthermore, there is a set $E\subset Q$ with measure $|E|=0$ such
that $u_k(x,t)\to u(x,t)$ if $(x,t)\in Q\setminus E$ so we
conclude that $u(x,t)=u^*(x,t)$ a.e. in $Q\cap Z^{-1}(\Omega)$.
This shows that $u(x,t)$ can be extended to a solution defined on
$Q\cup Z^{-1}(\Omega)$. Since $u^*(x,t)$ is smooth on
$Z^{-1}(\Omega)$ we will as usual modify $u(x,y)$ in a null set to
obtain $u(x,t)=u^*(x,t)$ everywhere on $Z^{-1}(\Omega)\cap Q$,
which means that we are picking a representative in the class of
$u\in L^\infty$ which restricted to $Z^{-1}(\Omega)\cap Q$ is
continuous. Assume without loss of generality that $m(0)\le0<M(0)$
and let's look at the boundary behavior of the holomorphic
function $U$ on $\Omega$. Since $U(\zeta)$ is bounded and $\Omega$
has a Lipschitz boundary the nontangential limit at a boundary
point $\zeta_0=x_0+i\varphi(x_0,0)$, $-\delta<x_0<\delta$, $$
\lim_{\Gamma(\zeta_0)\owns\zeta\to\zeta_0}
U(\zeta)=bU(\zeta_0)\in\ce\quad\text{ exists for a.e. }
\zeta_0\in\partial\Omega $$ where $\Gamma(\zeta_0)$ is a
nontangential region with vertex at $\zeta_0$. At points $\zeta_0$
where the limit does not exist we define $bU(\zeta_0)=0$ so $bU$
is now everywhere defined on $\partial\Omega$. We denote by
$\widehat U(\zeta)$ the natural extension to $\overline\Omega$ of
$U(\zeta)$, i.e., $$ \widehat U(\zeta)=\cases
U(\zeta),&\zeta\in\Omega,\\
                 bU(\zeta),&\zeta\in\partial\Omega, \endcases \tag1.3
$$
which is measurable and bounded and set $\widehat u(x,y)=\widehat U(x+i\varphi(x,y))$ for $(x,y)\in Q$. The proof of the theorem will rely heavily on the following representation formula for $u(x,y)$.
\enddemo
\proclaim{Lemma 1.3} Let $u(x,y)$, $U(\zeta)$ and $\widehat U(\zeta)$ be defined as above. After modification of $u(x,y)$  on a null set, the following identity holds
$$
u(x,y)=\widehat U(x+i\varphi(x,y))\qquad (x,y)\in Q.\tag 1.4
$$
\endproclaim
\demo{Proof}
We wish to prove that $u(x,y)=\widehat u(x,y)$ as distributions on $Q$. Let $(x_0,y_0)$ be an arbitrary point in $Q$. If $x_0+i\varphi(x_0,y_0)\in\Omega$ then $u(x,y)=
U\circ Z(x,y)=\widehat U\circ Z(x,y)=\widehat u(x,y)$ on a neighborhood of $(x_0,y_0)$. This shows that $u(x,y)=\widehat u(x,y)$ everywhere on $Q\cap Z^{-1}(\Omega)$. If $x_0+i\varphi(x_0,y_0)\in\partial\Omega$, then either $\varphi(x_0,y_0)=m(x_0)$ or $\varphi(x_0,y_0)=M(x_0)$. Let's assume first  that $m(x_0)=\varphi(x_0,y_0)<M(x_0)$. We consider the maximal vertical interval $\{x_0\}\times I(x_0)\subset Q$ that contains $(x_0,y_0)$ on which $\varphi(x,y)=m(x_0)$ and distinguish two cases.
\enddemo

\subheading{\bf Case 1} $I(x_0)=[y_0-\eta,y_0+\rho]$. Here
$\rho,\eta\ge0$ and $[y_0-\eta,y_0+\rho]\subset(0,\delta)$. Thus,
$\{x_0\}\times I(x_0)$ is contained in a rectangle
$R=(x_0-\mu,x_0+\mu)\times(y_0-\eta',y_0+\rho') \subset Q$,
$\eta'>\eta$, $\rho'>\rho$, such that $$
Z\Big(\big([x_0-\mu,x_0+\mu]\times\{y_0-\eta'\}\big)\cup\big(
[x_0-\mu,x_0+\mu]\times\{y_0+\rho'\} \big)\Big) \subset\Omega $$
and $$ \varphi(x,y)<M(x)\quad\text{on }\overline R. $$ In
particular, $\widehat u(x,y)=u(x,y)$ on a neighborhood of the
horizontal edges of the boundary of $R$. By the basic uniqueness
result on the Cauchy problem for locally integrable vector fields
mentioned at the beginning of this section, in order to  conclude
that $\widehat u=u$ as distributions on a neighborhood of
$\{x_0\}\times I(x_0)$ it will be enough to show that $L\widehat
u=0$ on $R$. For small $\epsilon>0$ consider the function $$
u_\epsilon(x,y)=U(x+i\varphi(x,y)+i\epsilon),\qquad (x,y)\in R. $$
While $U\circ Z$ is only defined on $Z^{-1}(\Omega)$, so  it fails
to be defined at points $(x,y)\in R$ such that
$\varphi(x,y)=m(x)$, this is not the case for $u_\epsilon(x,y)$
because the strict inequality $\varphi(x,y)+\epsilon>m(x)$  always
holds. If $(x,y)\in Q\cap Z^{-1}(\Omega)$ it is clear that
$\lim_{\epsilon\to0} u_\epsilon(x,y)=U\circ Z(x,y)=u(x,y)=
\widehat u(x,y)$ by the continuity of $U$. Let's study the limit
at a point $(x,y)\in R$, for which $\varphi(x,y)=m(x)$. Taking
account of (1.3) we see that $$ \lim_{\epsilon\searrow
0}u_\epsilon(x,y)=\lim_{\epsilon\searrow
0}U(x+i\varphi(x,0)+i\epsilon)=bU(x+i\varphi(x,0))=\widehat u(x,y)
$$ unless $x$ belongs to an exceptional set $E_1$ of measure
$|E_1|=0$. Since $Lu_\epsilon=0$, $u_\epsilon(x,y)\to \widehat
u(x,y)$ a.e. and the $u_{\epsilon}$ are uniformly bounded, we
conclude that $L\widehat u=0$ on $R$ which implies that $\widehat
u=u$ on a neighborhood of $(x_0,y_0)$.

\subheading{\bf Case 2} $I(x_0)=(0,y_0+\rho]$.

Here $\rho\ge0$ and $(0,y_0+\rho]\subset(0,\delta)$. Now
$\{x_0\}\times I(x_0)$ is contained in a rectangle
$R=(x_0-\mu,x_0+\mu)\times(0,y_0+\rho') \subset Q$, $\rho'>\rho$,
such that $$ Z\big(
[x_0-\mu,x_0+\mu]\times\{y_0+\rho'\}\big)\subset\Omega $$ and $$
\varphi(x,y)<M(x)\quad\text{on }\overline R. $$ This time we
conclude immediately that $\widehat u(x,y)=u(x,y)$ on a
neighborhood of the upper horizontal edge of the boundary of $R$.
But this is enough to repeat the argument of Case 1 and obtain via
uniqueness in the Cauchy problem that $\widehat u(x,y)=u(x,y)$ on
a neighborhood of $(x_0,y_0)$.

Note that $I(x_0)$ cannot be equal to $(0,\delta)$ because this
would imply that $m(x_0)=M(x_0)$. Hence, we have proved that if
$m(x_0)=\varphi(x_0,y_0)$ then $\widehat u=u$ on a neighborhood of
$(x_0,y_0)$. Similarly, we could prove that the same holds if
$\varphi(x_0,y_0)=M(x_0)$ by an analogous reasoning. This proves
the representation formula (1.4). \qed

We continue with the proof of the theorem assuming that $u(x,y)$
has been modified in a null set so (1.4) holds everywhere. As a
consequence  $u(x,y)$ is constant on the fibers
$F(x_0,y_0)=\{(x,y)\in Q:\quad
x=x_0,\,\,\varphi(x_0,y)=\varphi(x_0,y_0)\}$. In particular, if
$\varphi(x,y)$ is constant on some vertical segment of the form
$\{x_0\}\times(0,y_0]\subset Q$, the function $u(x,y)$  will also
be constant on that segment and it follows that $$
\lim_{\epsilon\searrow0}u(x_0,\epsilon)=u(x_0,y_0)=\widehat
u(x_0,y_0)=\widehat u(x_0,0) =\widehat U(x_0+i\varphi(x_0,0)) \tag
1.5 $$ We wish to see that (1.5) holds for almost all  points
$x_0\in(-\delta, \delta)$ such that $\varphi(x_0,0)=m(x_0)$. Let's
assume then that $m(x_0)=\varphi(x_0,0)<M(x_0)$. If
$m(x_0)=\varphi(x_0,y)$ for all $0<y\le y_0$ for some
$0<y_0<\delta$ we already saw the validity of (1.5) so we may
assume that there is a sequence $y_n\searrow0$ such that
$\varphi(x_0,0)<\varphi(x_0,y_n)<M(x_0)$ which implies that
$x_0+i\varphi(x_0,y_n)\in\Omega$. In this case we have $$
\lim_{y_n\searrow0}u(x_0,y_n)=\lim_{y_n\searrow0}U(x_0+i\varphi(x_0,y_n))=
bU(x_0+i\varphi(x_0,0)=\widehat U(x_0+i\varphi(x_0,0)) $$ unless
$x_0$ belongs to the exceptional set $E_1$ of measure $|E_1|=0$
introduced in the proof of Lemma 1.3.

Shrinking $\delta$ once again we may assume that
$\varphi(x,0)<M(x)$ for $|x|<\delta$. If
$m(x_0)<\varphi(x_0,0)<M(x_0)$ and $\rho$ is small enough,  the
extension of $u$ (still denoted by $u$) to $Q\cup Z^{-1}(\Omega)$
is defined for $|x-x_0|<\rho$, $-\rho<y<\rho$, for some small
$\rho$ such that $(x_0-\rho,x_0+\rho)\subset(-\delta,\delta)$, and
given by $u(x,y)=U(x+i\varphi(x,y))$. In particular, $u(x,y)$ is
smooth in a neighborhood of $(x_0,0)$ and  (1.5) follows by
continuity and the fact that $\widehat U(x_0+i\varphi(x_0,0))
=U(x_0+i\varphi(x_0,0))$ because $x_0+i\varphi(x_0,0)\in\Omega$.
Summing up, we have proved that $$
\lim_{\epsilon\searrow0}u(x,\epsilon)=\widehat U(x+i\varphi(x,0))
\qquad \text{a.e. } |x|<\delta. \tag 1.6 $$ Let
$\psi(x)\in\ccinf(-\delta,\delta)$. The dominated convergence
theorem gives $$ \<bu,\psi\>=\lim_{\epsilon\searrow0}\int
u(x,\epsilon)\psi(x)\,dx= \int \widehat
U(x+i\varphi(x,0))\psi(x)\,dx, $$ showing that $$ bu(x)=\widehat
U(x+i\varphi(x,0))\quad \text{a.e. }x\in(-\delta,\delta).\tag 1.7
$$ Consider now the domain $$
\Omega^+=\left\{\zeta=\xi+i\eta:\quad|\xi|<\delta,\,\,
\varphi(\xi,0)<\eta<M(\xi) \right\}, $$ and denote by $\widehat
U^+$  the restriction of $\widehat U$ to $\Omega^+$. Then (1.7)
may be rephrased as $$ bu(x)=bU^+(x+i\varphi(x,0))\quad \text{a.e.
}x\in(-\delta,\delta). $$ Indeed, it is clear that
$bU^+(\zeta)=bU(\zeta)$ if
$\zeta\in\partial\Omega\cap\partial\Omega^+$ and
$bU^+(\zeta)=U(\zeta)$ if $\zeta\in \Omega\cap\partial\Omega^+$.
We now invoke hypothesis (4) made on $bu(x)$. A change of
variables shows that $$
\int_{-\delta}^\delta\ln|bU^+(x)|\,dx=-\infty $$ and then it is
classical that $U^+$ must vanish identically, forcing $U$ and
therefore $u$ to vanish on $Q$.

So far we have proved  the first half  of the theorem. The second part will be given in the next section.

 \heading {2. End of the proof of Theorem 1.2}
 \endheading

We continue to work in the local coordinates $(x,y)$  that were used for the proof of the first part,  so the first integral has the form  $Z(x,y)=x+i\varphi(x,y)$ and
$$
L=X+iY=\frac{\partial}{\partial y}-\frac{i\varphi_y}{1+i\varphi_x}\frac{\partial}{\partial x},\qquad
$$
$$
X=\frac{\partial}{\partial y}-\frac{\varphi_y \varphi_x}{1+\varphi_x^2}\frac{\partial}{\partial x}, \qquad
Y=-\frac{\varphi_y }{1+\varphi_x^2}\frac{\partial}{\partial x}.
$$
It is easy to see that the hypothesis implies that $\varphi(0,y)=0$ for all $0\le y\le\delta$ for some $\delta>0$. Consider the function defined on $(-\delta,\delta)\times[0,\delta)$ by
$$
u(x,y)=\cases 0,&\text{for $x\le0$}\\
              \exp( Z^{-1}(x,y)),&\text{for $x>0$}.\endcases
$$ Since $Z(x,y)\not=0$ for $x\not=0$ it is clear that $u(x,y)$ is
smooth for $x\not=0$. Furthermore, as $x\searrow0$ we easily see
that $$
|u(x,y)|=\exp\left(-\frac{x}{x^2+\varphi^2(x,y)}\right)\longrightarrow
0 $$ because $\varphi(x,y)=O(x^2)$ for any  $0\le y<\delta$. The
same conclusion is valid for any derivative  $D_x^j D_y^k u(x,y)$
in place of $u(x,y)$. This proves that $u\in
C^\infty(-\delta,\delta)\times[0,\delta)$. Furthermore, $Lu=0$ for
$x\not=0$ because $u$ is constant for $x<0$ whereas  $u$ is a
holomorphic function of $Z$ defined on a neighborhood of
$Z((0,\delta)\times(0,\delta))$ for $x>0$. Since $Lu(x,y)$ is
continuous, it follows  that $Lu=0$ on
$(-\delta,\delta)\times(0,\delta)$. Also, since $u$ is continuous
up to $y=0$ it is apparent that $bu(x)=u(x,0)$. In particular,
$bu(x,0)=0$ for $x<0$ and $$ \int_{-\epsilon}^\epsilon
\ln|bu(x)|\,dx=-\infty $$ for any $\epsilon>0$. Finally,
$u(x,y)\not=0$ for any $x>0$. \qed

\heading{3. Pointwise convergence to the initial data}
\endheading
Always using the special coordinates $(x,y)$ defined in a
neighborhood of $[-a,a]\times[-b,b]$ where the first integral is
given by $Z=x+i\varphi(x,y)$, consider a function $u(x,y)$ defined
on $(-a,a)\times(0,b)$ that is measurable, bounded and satisfies
the equation $Lu=0$. Keeping the notation of Section 1, we set $$
M(x)=\sup_{0\le y\le b}\varphi(x,y),\qquad m(x)=\inf_{0\le y\le
b}\varphi(x,y), $$ and define $$ F=\{x\in(-a,a):\quad m(x)=M(x)\}.
$$ Then $F$ is closed and its complement  may be written as a
union of open intervals $$ (-a,a)\setminus F=\bigcup_{j=1}^\infty
(\alpha_j,\beta_j). $$ For any $j=1,2,\dots$ the set
$(\alpha_j,\beta_j)\times(0,b)$ is mapped by $Z$ into the domain
$$ \Omega_j=\left\{\zeta=\xi+i\eta:\quad \alpha_j<\xi<\beta_j,\,\,
m(\xi)<\eta<M(\xi) \right\}. \tag 3.1 $$ Since $m(x)<M(x)$ for
$\alpha_j<x<\beta_j$ the proof of Theorem 1.2 shows that there is
a holomorphic function $U_j$ defined on $\Omega_j$ having an
extension $\widehat U_j$ to $\overline\Omega$ such that
$u(x,y)=\widehat U_j\circ Z(x,y)$ for $\alpha_j<x<\beta_j$.
Furthermore, (1.6) holds, i.e., $$
\lim_{\epsilon\searrow0}u(x,\epsilon)=\widehat
U_j(x+i\varphi(x,0)) \qquad \text{a.e. } \alpha_j<x<\beta_j. $$ In
particular, the limit $\lim_{\epsilon\searrow0}u(x,\epsilon)$
exists for almost every $ x\in (-a,a)\setminus F$. On the other
hand, it is easy to see that the same limit exists for almost
every $x\in F$. Since this is a local property, it will be enough
to prove it in a neighborhood of a given point $x_0\in F$. To
simplify the notation let's assume that $x_0=0$. If $x\in F$,
$\varphi(x,0)=\varphi(x,y)$ for any $0<y<b$, so $Z(x,y)$ is
constant on $\{x\}\times(0,b)$, $x\in F$, and so is
$u_k(x,y)=P_k(Z(x,y))$, where $P_k(\zeta)$ is the sequence of
polynomials obtained from the Baouendi-Treves approximation scheme
used at the beginning of the proof of Theorem 1.2. Since
$u_k(x,y)\to u(x,y)$ a.e., for $|x|<\delta$, $0<y<\delta$ for some
$\delta>0$, there is a set $G$ of 2-dimensional measure zero such
that $u_k(x,y)\to u(x,y)$ for $(x,y)\notin G$. By Fubini's theorem
the 1-dimensional measure of $G_x=(\{x\}\times(0,\delta))\cap G$
is zero for a.e. $|x|<\delta$, i.e., for $x\notin
H\subset(-\delta,\delta)$, $|G_x|=0$, with $|H|=0$. If $x\in
F\setminus H$, $u_k(x,y)\to u(x,y)$ for a.e.~$y$ and we see that
$y\mapsto u(x,y)$ is constant almost everywhere. Thus, modifying
$u(x,y)$ in a set of 2-dimensional measure zero, we obtain that
$y\to u(x,y)$ is constant for $x\in F\setminus H$ and
$\lim_{\epsilon\to0}u(x,\epsilon)$ exist for $x\in F\setminus H$.
Observe that this limit has to equal $bu(x)$ a.e. since we can
find a sequence $y_j\mapsto 0$ such that the traces $u(.,y_j)$
converge to $bu$ weakly in $L^{\infty}(-a,a)$. Hence, we have
proved the existence of the limit
$\lim_{\epsilon\to0}u(x,\epsilon)$ for a.e. $x\in(-\delta,\delta)$
and therefore for a.e. $x\in(-a,a)$. In particular, this gives an
alternative proof of Lemma 1.1 when $L$ is locally integrable.

In the convergence result just proved, the boundary point $(x_0,0)$ is approached  perpendicularly to initial curve $\{y=0\}$ (in the special coordinates we are using). Since in the classical case where $L$ is the Cauchy-Riemann operator and $u$ is holomorphic the limit is still valid on nontangential regions of approach, it is natural to try to replace the normal set of approach $\{(x_0,\epsilon)\}$ by a larger set. Keeping in mind that $L$ is elliptic exactly at the points where $\varphi_y\neq0$, we di
tinguish two types of points $x_0\in(-a,a)$:
\roster
\item"(I)" there exists $\eta=\eta(x_0)>0$ such that $\varphi_y(x_0,y)=0$ for $0\le y\le\eta$;
\item"(II)" there exists a sequence of points $y_n>0$ such that $y_n\searrow0$ and $\varphi_y(x_0,y_n)\neq 0$.
\endroster
We will attach to every point $x_0$ a set of approach
$\Gamma(x_0)$. In case (I) we keep the normal set of approach and
simply define $\Gamma(x_0)=\{x_0\}\times(0,b)$ . In case (II) we
proceed as follows: since the hypothesis implies that
$m(x_0)<M(x_0)$ it follows that $x_0\in(\alpha_j,\beta_j)$ for
some $j$ and $u(x,y)=\widehat U_j\circ Z(x,y)$ for
$\alpha_j<x<\beta_j$ with $U_j$ holomorphic and bounded in the
open set $\Omega_j$ given by (3.1). Since $U(\zeta)\to
bu(x_0+i\varphi(x_0,0))$ nontangentially unless $x_0$ belongs to
an exceptional set of measure zero, we define $\Gamma(x_0)$ as the
interior of $$ Z^{-1}\big(\{\xi+i\eta:\quad
|\xi-x_0|\le|\eta-\varphi(x_0,0)|\}\big)\cap\{y>0\} $$ which is an
open set that contains $\{x_0\}\times(0,b)$ and for almost every
$x_0\in(\alpha_j,\beta_j)$ satisfies $$
\lim_{\Gamma(x_0)\owns(x,y)\to(x_0,0)}u(x,y)=bu(x_0), $$ as follow
from the  convergence  of $U_j(\zeta)\to\zeta_0=x_0 +
i\varphi(x_0,0)$ as $\zeta=\xi+i\eta\to x_0 + i\varphi(x_0,0)$,
when $|\eta-\varphi(x_0,0)|\le |\xi-x_0|$. Note that if $L$ is
elliptic at $(x_0,0)$, i.e., $\varphi_y(x_0,0)\neq0$, then
$\Gamma(x_0)$ contains a cone $y>\mu|x-x_0|$, with $\mu>0$, so we
recover the classical nontangential convergence valid for elliptic
vector fields. On the other hand, if $\varphi_y(x_0,0)=0$ and
$x_0$ is of type (II)) then $\Gamma(x_0)$ is still an open
neighborhood of $\{x_0\}\times(0,b)$ which cannot contain any cone
$y>\mu|x-x_0|$, $\mu>0$, because its width  at height $y$ is
$O(|y|^2)$ and it is contained in a cuspoidal region
$y>\mu|x-x_0|^{1/2}$ with vertex $x_0$. We summarize these facts
in a more invariant way as

\proclaim{Theorem 3.1} Let $$ L=\frac{\partial}{\partial
y}+a(x,y)\frac{\partial}{\partial x}=X+iY, $$ $a(x,t)\in
C^{\infty}$ on $(-a,a)\times(-b,b)$, be a locally integrable
vector field.  Denote by $\gamma_x$ the integral curve of $X=\re
L$ that stems from $\{x\}$ and enters $\{y>0\}$. Then to each
$(x,0)$ we may associate a set of convergence $\Gamma(x)$, which
is an open neighborhood of $\gamma_x$ if $(x,0)$ is of type (II)
and reduces to $\gamma_x$ if $(x,0)$ is of type (I) so that for
any  $u(x,y)\in L^\infty((-a,a)\times(0,b))$ that satisfies
$Lu(x,y)=0$ for $y>0$ we have $$
\lim_{\Gamma(x_0)\owns(x,y)\to(x_0,0)}u(x,y)=bu(x_0),\qquad
\text{\rm a.e. }x_0\in(-a,a). $$
\endproclaim

In the example below there  are points $(x_0,0)$ of type (I) and type (II). The points of type
(II) are all elliptic points and the solution $u(x,y)$ converges nontangentially to its
 boundary values $bu(x)$ at those points while the convergence occurs strictly along $\gamma_x$
 at points $(x,0)$ of type (I). The example shows that in general we cannot hope to enlarge
 $\Gamma(x)$ to an open  neighborhood of $\gamma_x$ at points of type (I). We refer the reader to
 [BH2] for related results on pointwise convergence to the boundary value.

\example{Example}
Let $K\subset(-1,1)\subset\erre$ be a Cantor set with positive measure $|K|>0$ and denote by
$(\alpha_j,\beta_j)$, $j=1,2,\dots$, its complementary intervals in $(-1,1)$.
 Let $b(x)\in\cinf(-1,1)$ such that
$$
\cases
b(x)=0&\text{if $x\in K$ and}\\
b(x)>0,&\text{if $x\in (-1,1)\setminus K$}.
\endcases %\tag 2.1
$$
Define
$$
Z(x,y)=x+iyb(x);\qquad L=\frac{\partial}{\partial y}-\frac{ib(x)} {1+iyb'(x)}
\frac{\partial}{\partial x},
$$
where $b'(x)$ denotes the derivative of $b(x)$. It is readily checked that $LZ=0$ and $Z_x$ does not vanish so $L$ is locally integrable. Note that $L$ is nonelliptic exactly at the points $(x,y)\in(-1,1)\times\erre$ where $b(x)=0$, i.e., at the points of $K\times\erre$. Thus, the points of $K$ are of type (I) and those of $(-1,1)\setminus K$ are elliptic points of type (II).

Consider the characteristic function of the set $K$, $\chi(x)=1$
if $x\in K$ and $\chi(x)=0$ otherwise and set $u(x,y)=\chi(x)$.
Clearly, $u\in L^\infty(\erre^2)$ and since $u(x,y)$ is
independent of $y$ it follows that $\partial_y u=0$. Furthermore,
$\chi'(x)$ is a distribution supported in $K$ and since $b(x)$
vanishes to infinite order on $K$ it follows that $b\chi'=0$. This
shows that $Lu=0$. Finally, we point out that if $x_0\in K$, any
neighborhood of $\{x_0\}\times(0,1)$ in $\{y>0\}$ contains a
sequence $(x_n,y_n)$ of points converging to $(x_0,0)$ such that
$u(x_n,y_n)=0$ while $bu(x_0)=1$.

\endexample

\heading{ Appendix}
\endheading

\demo{Proof of Lemma {\sl 1.1}}  Since $f$ is  bounded, we need
only show that  $f(x,y)$ converges in ${\D}'(-a,a)$ to a
distribution $bf(x)$ as $y\searrow 0$ as the weak compactness of
the unit ball of $L^\infty(-a,a)$ will then show that $bf\in
L^\infty(-a,a)$. We will proceed as in [BH1] with minor
modifications. Let $\phi \in C_0^{\infty}(-a,a)$. For $\epsilon
\geq 0$ sufficiently small, set $$
L^{\epsilon}=\frac{\partial}{\partial
y}+a(x,y+\epsilon)\frac{\partial}{\partial x} $$ We will choose
$\phi_0^{\epsilon}$ and $\phi_1^{\epsilon} \in
C^{\infty}((-a,a)\times[0,b))$ such that if $$
\Phi^{\epsilon}(x,y)=\phi_0^{\epsilon}(x,y)+y\phi_1^{\epsilon}(x,y),
$$ then $$ (1)\quad \Phi^{\epsilon}(x,0)=\phi(x),\quad \text{and }
\quad (2)\quad |(L^{\epsilon})^{\ast}\Phi^{\epsilon}(x,y)| \leq
Cy, $$ where $(L^{\epsilon})^{\ast}$ denotes the formal transpose
of $L^{\epsilon}$ and $C$ depends only on the derivatives of
$\phi$ up to order $2$. In particular, $C$ will be independent of
$\epsilon$. Define $\phi_0^{\epsilon}(x,y)=\phi(x)$ and  write $$
L^{\epsilon}=\frac{\partial}{\partial
y}+Q^{\epsilon}(x,y,\frac{\partial}{\partial x}) , $$ and define
$$ \phi_1^{\epsilon}(x,y)=-\frac{\partial}{\partial y}
\phi_{0}^{\epsilon}(x,y) +(Q^{\epsilon})^{\ast}\phi_{0}^{\epsilon}
$$ One easily checks that (1) and (2) above hold with these
choices of the $\phi_j^{\epsilon}$. We will next use the
integration by parts formula of the form $$ \int
u(x,T)w(x,T)\,dx-\int u(x,0)w(x,0)\,dx=\int_0^T\int_{{\Bbb
R}^n}(wPu-uP^{\ast}w)\,dxdy $$ which is valid for $P$ a vector
field, $u$ and $w$ in $C^1({\Bbb R}\times [0,T])$ and the
$x-$support of $w$ contained in a compact set in ${\Bbb R}$. Note
that the $x$-support of $\Phi^{\epsilon}(x,y)$ is contained in the
support of $\phi(x)$. Let $\psi \in C_0^{\infty}(B_1(0))$, where
$B_1(0)$ denotes the ball of radius $1$ centered at the origin in
${\Bbb R}^2$, be of the form $\psi(x,y)=\alpha(x)\beta(y)$. Assume
$\int \alpha(x)\, dx=\int \beta(y)\,dy=1$, and for $\delta >0$,
let $\psi_{\delta}(x,y)=\delta
^{-2}\psi(x/\delta),y/\delta)=\alpha_\delta(x)\beta_\delta(y)$.
For $\epsilon>0$, set $f_{\epsilon}(x,y)=f(x,y+\epsilon)$. Observe
that if $\delta < \epsilon$, then the convolution
$f_{\epsilon}*\psi_{\delta}(x,y)$ is $C^{\infty}$ in the region
$y\geq 0$. In the integration by parts formula above set
$u(x,y)=f_{\epsilon}*\psi_{\delta}(x,y)$, $w(x,y)=
\Phi^{\epsilon}(x,y)$ and $P=L^{\epsilon}$. We get $$ \align
  \int_{-a}^af_{\epsilon}*\psi_{\delta}(x,0) \phi(x)\,dx &=\int_{-a}^a f_{\epsilon}*\psi_{\delta}(x,T) \Phi^{\epsilon}(x,T)\,dx \\
                               &\qquad-\int_0^T\int_{-a}^aL^{\epsilon}\left (f_{\epsilon}*\psi_{\delta}\right )\Phi^{\epsilon}\,dxdy \tag a.1 \\
                               &\qquad\qquad+\int_0^T\int_{-a}^a f_{\epsilon}*\psi_{\delta}(L^{\epsilon})^{\ast}\Phi^{\epsilon}\,dxdy.
  \endalign
$$ We have chosen $0<T<b$ such that $x\mapsto f(x,T)$ is bounded
and measurable. Fix $\epsilon >0$. Let $\delta \rightarrow 0^+$.
Note that $\{f_{\epsilon}*\psi_{\delta}(x,y)\}$ is uniformly
bounded and converges almost everywhere  to $f_{\epsilon}(x,y)$ on
a neighborhood $W$ of $\text{ supp } \phi \times [0,T]$. Hence, $$
L^{\epsilon}\left(f_{\epsilon}*\psi_{\delta}\right )\rightarrow
L^{\epsilon}f_{\epsilon} $$ in ${\D}'(W)$ as $\delta \rightarrow
0^+$. Moreover, $L^{\epsilon}f_{\epsilon}(x,y)=Lf(x,y+\epsilon)\in
L^{\infty}$. Hence, by Friederichs' Lemma, $$
L^{\epsilon}\left(f_{\epsilon}*\psi_{\delta}\right )\rightarrow
L^{\epsilon}f_{\epsilon} $$ in $L^2(W)$ as $\delta \rightarrow
0^+$. Furthermore, using that the trace map $(0,b)\owns y\mapsto
f(x,y)\in\D'((-a,a))$ is continuous, we also see that the first
integral on the right hand side of (a.1) converges to $$
\int_{-a}^a f(x,T+\epsilon) \Phi^{\epsilon}(x,T)\,dx. $$ We thus
get $$ \align
  \int_{-a}^af(x,\epsilon) \phi(x)\,dx &=\int_{-a}^a f(x,T+\epsilon) \Phi^{\epsilon}(x,T)
       \,dx \\
                               &\qquad -\int_0^T\int_{-a}^aL^{\epsilon}f_{\epsilon}(x,y)                        \Phi^{\epsilon}(x,y)\,dxdy  \\
                               &\qquad\qquad+\int_0^T\int_{-a}^a f_{\epsilon}(x,y)(L^{\epsilon})^{\ast}\Phi^{\epsilon}(x,y)\,dxdy
  \endalign
$$
In the third integral on the right, we have
 $$
|f_{\epsilon}(x,y)(L^{\epsilon})^{\ast}\Phi^{\epsilon}(x,y)|\leq Cy ,
$$
 where $C$ depends only on the derivatives of $\phi$ upto order $2$.
 By the dominated convergence theorem, as $\epsilon \to 0$,
 this third integral converges to
$\int_0^T\int_{-a}^a fL^{\ast}\Phi^0\,dxdy$. In the second integral on the right,
note that since $Lf\in L^2(X\times (0,T)) $, as $\epsilon\to 0$, the
translates $L^{\epsilon}f_{\epsilon}=(Lf)_{\epsilon}\to Lf$ in $L^2$. We thus get
$$
\langle
bf,\phi \rangle = \int_{-a}^a f(x,T)\Phi(x,T)\,ds-\int_0^T\int_{-a}^a Lf\Phi\, dx dy+
\int_0^T\int_{-a}^a fL^{\ast}\Phi\,dxdy ,
$$
where $\Phi=\Phi^0$.
\qed
 \enddemo


\Refs
\widestnumber\key{BH2}
\refstyle{A}

\ref
\key BT
\by M. S. Baouendi and F. Treves
\paper A property of the functions and distributions annihilated by a locally integrable system of complex vector fields
\jour Ann. of Math.
\vol 113
\yr 1981
\pages 387--421
\endref

\ref
\key BH1
\by S. Berhanu and J.Hounie
\paper  An F. and M. Riesz theorem for planar vector fields
\vol 320
\yr 2001
\pages 463-485
\jour Math. Ann.
\endref

\ref
\key BH2
\by S. Berhanu and J.Hounie
\paper  On boundary properties of solutions of complex vector fields
\jour preprint
\endref

\ref
\key C
\by  P. Cohen
\paper The non-uniqueness of the Cauchy problem
\jour O.N.R. Tech. Report
\vol 93
\moreref\jour Stanford
\yr  1960
%\pages 1324--1327
\endref

\ref
\key Co
\by P.Cordaro
\paper Approximate solutions in locally
 integrable structures
\jour Fields Institute Communications volume:
  Differential Equations and Dynamical Systems
  (in honor to Waldyr Muniz Oliva), to appear
\endref

\ref
\key Du
\by P. Duren
\book Theory of $H^p$ spaces
\publ
Academic Press
\yr 1970
\endref

\ref \key F \by P. Fatou \paper S\'eries trigonom\'etriques e
s\'eries de Taylor \jour Acta Math. \vol 30 \yr 1906 \pages
335--400
\endref

\ref
\key HM
\by J. Hounie and P. Malagutti
\paper On the convergence of the Baouendi-Treves approximation formula
\jour Comm. P.D.E.
\vol 23
\yr 1998
\pages 1305--1347
\endref

\ref
\key HT
\by J. Hounie and J. Tavares
\paper On removable singularities of locally solvable differential operators
\jour Invent. Math.
\vol 126
\yr 1996
\pages 589--623
\endref

\ref
\key Hor
\by L. H\"ormander
\book The Analysis of linear partial differential operators I
\publ
Springer-Verlag
\yr 1990
\endref

\ref
\key J
\by B.J\"oricke
\paper Deformation of CR manifolds, minimal points and CR manifolds
 with the microlocal analytic extension property
\jour J. Geom. Anal.
\yr 1996
\pages 555--611
\endref

\ref
\key RR
\by F. Riesz and M. Riesz
\paper \"Uber die Randwerte einer analytischen Funktion
\jour Quatri\`eme Congr\`es de Math. Scand. Stockholm
\yr 1916
\pages 27--44
\endref


\ref
\key T1
\by F. Treves
\book Hypo-analytic structures, local theory
\publ Princeton University Press
\yr 1992
\endref

\ref
\key  T2
\by F. Treves
\paper Approximation and representation of solutions in locally integrable
structures with boundary
\jour Aspects of Math. and Applications
\yr 1986
\pages 781-816
\endref

\ref
\key  T3
\by F. Treves
\book Approximation  and representation of functions and  distributions annihilated by a system of complex vector fields
\publ Centre de Math\'ematiques.\'Ecole   Polytechnique,
Pa\-lai\-seau, France
\yr 1981
\endref

\ref
\key Z
\by C. Zuily
\book Uniqueness and non-uniqueness in the Cauchy Problem
\publ Birkh\"auser, Boston-Basel-Stuttgart
\yr 1983
\endref


\endRefs
\enddocument